\documentclass[a4paper,11pt]{amsart}
\usepackage[applemac]{inputenc}

\newcommand{\R}{\mathbb{R}}
\newcommand{\Z}{\mathbb{Z}}
\newcommand{\N}{\mathbb{N}}
\newtheorem{theorem}{Theorem}
\newtheorem{proposition}{Proposition}

\theoremstyle{remark}
\newtheorem{remark}{Remark}

\def\div{\, \mbox{div}\,  }    \def \e {\varepsilon} 
\def\tilde{\widetilde} 

\def\dq{\Delta_q} 
 
\def\d{\partial}

\def\cC{{\mathcal C}}

\def\div{\, \mbox{\rm div}\,  }

\newenvironment{p}{ 
\noindent\textit{\textbf{Proof:}}~} 
{\hfill\rule{2.1mm}{2.1mm} 
\smallbreak} 
\newcommand{\Int}{\displaystyle \int} 
\newcommand{\Frac}{\displaystyle \frac}

\newcommand{\du}{\delta\! u} 
\newcommand{\dv}{\delta\! v} 
 
\newcommand{\dt}{\delta\!\theta} 
 
\newcommand\dPi{\delta\!\Pi} 
 
\def \epsilon {\varepsilon}        
\begin{document}
\title{ Global well-posedness issues for the inviscid Boussinesq system with 
Yudovich's type data}
\author{Rapha\"el Danchin$^1$ and Marius Paicu$^2$}
\thanks{$^1$Universit\'e Paris-Est, Laboratoire d'Analyse 
et de Math\'ematiques Appliqu\'ees, UMR 8050,
 61 avenue du G\'en\'eral de Gaulle,
94010 Cr\'eteil Cedex, France.  E-mail: danchin$@$univ-paris12.fr}
\thanks{$^2$Universit\'e Paris 11, Laboratoire de Math\'ematiques,  
 B\^atiment 425, 91405 Orsay Cedex, France. 
 E-mail: marius.paicu$@$math.u-psud.fr}

\date\today
\begin{abstract}
The present paper is dedicated to the study of the global existence
for the inviscid two-dimensional Boussinesq system. 
We focus on finite energy data with bounded vorticity
and we find out that, under quite a natural additional assumption on 
the initial temperature, there exists a global unique solution. None smallness conditions are imposed on the data.
The global existence issues for  infinite energy initial velocity,
and for the B\'enard system are also discussed. 
\end{abstract}

\maketitle

 \section*{Introduction}

The incompressible  Euler equations 
have been intensively studied from a mathematical viewpoint. 
The present paper aims at extending   
the celebrated result by Yudovich concerning the two-dimensional 
Euler system (see \cite{Yudo}) to the 
 following \emph{two-dimensional Boussinesq system}:
$$
\begin{cases}
\partial_t\theta+u\cdot\nabla \theta-\kappa \Delta\theta=0\\
\partial_t u+u\cdot\nabla u-\nu\Delta u+\nabla\Pi=\theta\, e_2
\quad\hbox{with}\quad e_2=(0,1),\\
\div u=0.
\end{cases}\leqno(B_{\kappa,\nu})
$$
The above system describes the evolution of the velocity field
$u$ of a two-dimensional 
incompressible fluid moving under a vertical force the  magnitude $\theta$ 
of which is transported with or without diffusion by $u.$
Above the molecular diffusion parameter $\kappa$ and
viscosity $\nu$ are nonnegative, and $\Pi$ stands for the pressure 
in the fluid. 
For the sake of simplicity, we restrict our attention
to the case where the space variable $x$ belongs to 
the whole plan $\R^2$ (our results extend with no difficulty
to periodic boundary conditions, though).

The Boussinesq system is of relevance to study a number
of models coming from atmospheric or oceanographic turbulence
where rotation and stratification play an important role (see e.g. \cite{Ped}).
The scalar function $\theta$ may for instance represent temperature
variation in a gravity field, and $\theta\,e_2,$ the buoyancy force.

{}From the mathematical point of view, if both $\kappa$ and $\nu$ are
\emph{positive} then standard energy methods yield global existence of smooth
solutions for arbitrarily large data (see e.g. \cite{CD,Guo}). In contrast, in
the case  when $\kappa=\nu=0,$ the Boussinesq system exhibits vorticity
intensification and   the global well-posedness issue  remains an unsolved 
challenging open problem (except if $\theta_0$ is a constant of course)
which may be formally compared to the similar problem 
for the three-dimensional axisymmetric Euler equations \emph{with swirl}
(see e.g. \cite{ES} for more explanations). 

As pointed out by H. K. Moffatt in \cite{Moffatt}, knowing whether
having $\kappa>0$ or $\nu>0$ precludes the formation of finite time singularities
is an important issue. 
In \cite{DP1}, we stated that in the case $\kappa=0$
and $\nu>0$ no such formation may be encountered for finite energy initial data.
More precisely, we stated that for any $(\theta_0,u_0)$ in $L^2(\R^2)$ with $\div u_0=0,$
System $(B_{0,\nu})$ has a unique global finite energy solution.

In the present paper,  we aim at investigating the opposite case, namely 
$\kappa>0$ and $\nu=0.$
The corresponding Boussinesq system thus reads 
$$\begin{cases} \partial_t\theta+u\cdot\nabla
\theta-\kappa\Delta \theta=0\\ \partial_t u+u\cdot\nabla u+\nabla\Pi=\theta\,
e_2\\ \div u=0\end{cases}\leqno(B_{\kappa,0}) 
$$
and may be seen as a coupling between the two-dimensional Euler equations
and a transport-diffusion equation. 
In passing, let us point out that in the case $\theta\equiv0,$
 System $(B_{\kappa,0})$  reduces to the Euler equation. 

It is well known that the standard Euler equation is globally well-posed in
$H^s$ for any $s>2.$
A similar result has been stated for $(B_{\kappa,0})$
in the case $s\geq3$ by D. Chae in \cite{dingo}, then
 extended to rough data by T. Hmidi and S. Keraani in 
\cite{HK2}. There, global well-posedness is shown whenever the initial
velocity $u_0$ belongs to $B^{1+\frac2p}_{p,1}$ and 
the initial temperature $\theta_0$ is in $L^r$ for
some $(p,r)$ satisfying $2<r\leq p\leq\infty$ (plus a technical condition
if $p=r=\infty$). 
Let us emphasize that in the Besov spaces framework, the assumption 
on $u_0$ is somewhat optimal (since it is optimal for the
standard Euler equations, see \cite{Vishik}).

Here we want to state  global existence 
for less regular data satisfying Yudovich's type conditions. 
Roughly, we want to consider data $(\theta_0,u_0)$ in $L^2$ 
such that the initial vorticity $\omega_0:=\d_1v_0^2-\d_2v_0^1$ is bounded. 
Note however that, since we expect the corresponding solution to have
bounded vorticity for all positive time, we have to introduce an
additional assumption on $\theta_0.$
Indeed, the vorticity equation reads
$$
\d_t\omega+u\cdot\nabla\omega=\d_1\theta.
$$
Therefore, since no gain of smoothness 
may be expected from the above transport equation,
having $\omega$ bounded requires that 
$\d_1\theta\in L^1_{loc}(\R_+;L^\infty).$
Now, considering that $\theta$ satisfies the following heat equation
$$
\d_t\theta-\kappa\Delta\theta=f\quad\hbox{with}\quad f:=-u\cdot\nabla\theta,
$$
 the assumptions on $\theta_0$ should ensure
that 
\begin{equation}\label{eq:theta0}
\nabla e^{\kappa t\Delta}\theta_0\in L^1_{loc}(\R_+;L^\infty)
\end{equation}
where $(e^{\lambda\Delta})_{\lambda>0}$ stands  for the heat semi-group.
\smallbreak
It turns out that \eqref{eq:theta0} is equivalent to 
having $\nabla\theta_0$ in the \emph{nonhomogeneous Besov space}
$B^{-2}_{\infty,1}$ (see e.g. \cite{BCD}).
This motivates the following statement which is the main 
result of the paper:
\begin{theorem}\label{th:main1} Let  $\theta_0\in L^2\cap
B^{-1}_{\infty,1}$ and  $u_0\in L^2$  with $\div u_0=0.$
Assume in addition that the initial vorticity $\omega_0$
belongs to $L^r\cap L^\infty$ for some $r\geq2.$
 System  $(B_{\kappa,0})$ admits a unique
global solution  $(\theta,u)$ satisfying  
\begin{equation}\label{eq:regularity}
\begin{array}{c}\theta\in\cC(\R_+;L^2\cap
B^{-1}_{\infty,1}) \cap L^2_{loc}(\R_+;H^1)\cap
L^1_{loc}(\R_+;B^1_{\infty,1}),\\[1.5ex]
 u\in\cC^{0,1}_{loc}(\R_+;L^2) \ \hbox{
and }\ \omega\in L^\infty_{loc}(\R_+;L^r\cap
L^\infty).\end{array}\end{equation}
\end{theorem}  
\begin{remark} As a by-product of our proof, 
we gather  that if in addition  $\theta_0\in L^p$ 
(resp. $u_0\in B^1_{\infty,1}$)  for some
$p\in[1,+\infty]$ then $\theta\in L^\infty(\R_+;L^p)$
(resp. $u\in\cC(\R_+;B^1_{\infty,1})$). 
\end{remark}
\begin{remark} The $B^{-1}_{\infty,1}$
hypothesis over $\theta_0$ is quite  mild  compared to 
the $L^2$ hypothesis. Indeed, it may be shown that 
$L^2$ is continuously embedded in the Besov space $B^{-1}_{\infty,2}$
which is slightly larger than $B^{-1}_{\infty,1}.$
\end{remark}
The paper unfolds as follows. 
In the first section, we prove Theorem \ref{th:main1}. 
In the second section, motivated by the fact that having $u_0$ in $L^2$ and 
$\omega_0\in L^1$ requires the vorticity
to have zero average over $\R^2,$  we consider initial velocities
which are $L^2$ perturbations of infinite energy smooth stationary solutions for 
the incompressible Euler equations. 
Some extensions to Theorem \ref{th:main1} are discussed
in the third section. 
A few technical inequalities have been postponed in the appendix.


\section{Proof of Theorem \ref{th:main1}}

Proving Theorem \ref{th:main1} requires our using
the (nonhomogeneous) Littlewood-Paley decomposition.
One can proceed as in \cite{C}:
first we consider a dyadic partition of unity: 
$$1=\chi(\xi)+\sum\limits_{q\geq 0}\varphi(2^{-q}\xi),$$ 
for some nonnegative function 
 $\chi\in\cC^\infty(B(0,\frac 43))$ with value  $1$ over the ball
$B(0,\frac 34),$  and $\varphi(\xi):=\chi(\xi/2)-\chi(\xi)$.
\smallbreak
Next, we introduce the \emph{dyadic blocks}
$\dq$ of our decomposition by setting 
$$
\dq u:=0\ \hbox{ if }\ q\leq-2,\qquad
\Delta_{-1}u:={\mathcal F}^{-1}(\chi{\mathcal F}u)
\quad\hbox{and}\quad
\dq u:={\mathcal F}^{-1}(\varphi(2^{-q}\cdot){\mathcal F} u)\ \text{ if }
\  q\geq 0.
$$
One may prove  that for all tempered
distribution $u$ the following \emph{Littlewood-Paley decomposition} holds
true:  $$u=\sum\limits_{q\geq -1}\Delta_q u.$$ 
For $s\in\R,$ $p\in[1,\infty]$ and  $r\in[1,\infty],$ 
one can now define the \emph{nonhomogeneous Besov space} 
 $B^s_{p,r}:=B^s_{p,r}(\R^2)$ as the set of tempered distributions  $u$ over
$\R^2$ so that  $$\|u\|_{B^s_{p,r}(\R^2)}:=\bigl\|2^{qs}\|\Delta_q
u\|_{L^p(\R^2)}\bigr\|_{\ell^r(\Z)}<\infty.$$
We shall also use several times the following 
well-known fact for incompressible fluid mechanics (see  the proof in e.g.
\cite{C}, Chap. 3): 
\begin{proposition}\label{p:biot-savart} For any 
$p\in]1,\infty[$ the operator $\omega\mapsto\nabla u$
is bounded in $L^p.$ More precisely, there exists a constant $C$
such that 
$$
\|\nabla u\|_{L^p}\leq C\frac{p^2}{p-1}\,\|\omega\|_{L^p}.
$$
\end{proposition}
One can now tackle the proof of Theorem \ref{th:main1}.
One  shall proceed as follows.
\begin{enumerate}
\item[1.] We smooth out the data so as to get a sequence
of global  smooth solutions to $(B_{\kappa,0}).$
\item[2.] Energy estimates are proved.
\item[3.] We establish estimates in larger norms. 
\item[4.] We state uniform estimates for the first order time derivatives.
\item[5.] We pass to the limit in the system by means of compactness arguments.
\item[6.] Uniqueness is proved.
\end{enumerate}

\subsubsection*{First step.} We  smooth out the initial data
  $(\theta_0,u_0)$ (use e.g. a convolution process)
  and get a sequence  of smooth initial data
   $(\theta_0^n,u_0^n)_{n\in\N}$ which is bounded in the space
   given in the statement of the theorem. 
In addition, those smooth data belong to all the Sobolev spaces
$H^s.$ Hence, applying Chae's result \cite{dingo} provides us with 
a sequence of smooth global solutions  
 $(\theta^n,u^n)_{n\in\N}$
 which belong to all the spaces  $\cC(\R_+;H^s).$
 {}From system  $(B_{\kappa,0})$ and standard product laws in Sobolev spaces,
 we deduce that $(\theta^n,u^n)$ belongs to 
$\cC^1(\R_+;H^s)$  for all $s\in\R,$ and thus also 
to  $\cC^1(\R_+;L^p)$ for all $p\in[2,\infty].$
This will be more than enough to make the  computations in the following two
steps rigorous. 
\subsubsection*{Second step.} We want to state  energy type estimates
for $(\theta^n,u^n).$
Let us first take the $L^2(\R^2)$ inner product of $\theta^n$ with  the equation
satisfied by $\theta^n.$ Performing 
a space integration by parts in the diffusion term and a time integration yields
\begin{equation}\label{eq:thetaL2}
\|\theta^n(t)\|_{L^2}^2+2\kappa\int_0^t\|\nabla\theta^n\|_{L^2}^2\,d\tau
=\|\theta^n_0\|_{L^2}^2\quad\hbox{for all}\ t\in\R_+.
\end{equation}
As for the velocity $u^n,$ a similar  argument gives
$$
\|u^n(t)\|_{L^2}\leq\|u^n_0\|_{L^2}+\int_0^t\|\theta^n\|_{L^2}\,d\tau.
$$
Hence, bounding $\|\theta^n\|_{L^2}$ according to \eqref{eq:thetaL2}, we get
 \begin{equation}\label{eq:uL2}
\|u^n(t)\|_{L^2}\leq\|u_0^n\|_{L^2}+t\|\theta^n_0\|_{L^2}.
\end{equation}

\subsubsection*{Third step.}
This is the core of the proof of global existence. 
We here want to get uniform estimates for the Besov norms
of $\theta^n$ and for $\|\omega^n\|_{L^r\cap L^\infty}.$

Let us first consider the vorticity. As explained in the introduction, 
we have 
$$
\d_t\omega^n+u^n\cdot\nabla\omega^n=\d_1\theta^n.
$$
Therefore, for all  $p\in[r,\infty],$
\begin{equation}\label{eq:tourbillonLp}
\|\omega^n(t)\|_{L^p}\leq
\|\omega^n_0\|_{L^p}+\int_0^t\|\d_1\theta^n\|_{L^p}\,d\tau.
\end{equation} 
Hence,   getting uniform bounds on  $\|\omega^n\|_{L^r\cap L^\infty}$
requires  uniform bounds for 
 $\d_1\theta^n$ in  the space $L^1_{loc}(\R_+;L^r\cap L^\infty).$
 Because Equality \eqref{eq:thetaL2} supplies a  bound in $L^2(\R_+;L^2)$ 
for  $\d_1\theta^n,$  it is enough to get a suitable bound for
$(\d_1\theta^n)_{n\in\N}$  in  $L^1_{loc}(\R_+;L^\infty).$
 Given that the operator  $\d_1$ maps 
 $B^1_{\infty,1}$ in  $B^0_{\infty,1},$ and that 
 $B^0_{\infty,1}\hookrightarrow L^\infty,$ 
the problem reduces to proving uniform  estimates for 
$\theta^n$ in $L^1_{loc}(\R_+;B^1_{\infty,1}).$

For doing so, we rewrite the equation for $\theta^n$ as follows~:
\begin{equation}\label{eq:star}
\d_t\theta^n-\kappa\Delta\theta^n=-u^n\cdot\nabla\theta^n
\end{equation}
and take advantage of the smoothing properties of
the heat equation. 
More precisely, it is stated in the appendix that for all 
$\alpha\in[1,\infty],$ \begin{equation}\label{eq:chaleur}
\kappa^{\frac1\alpha}\|\theta^n\|_{L_T^\alpha(B^{-1+\frac2\alpha}_{\infty,1})}
\leq C(1+\kappa t)^{\frac1\alpha}\biggl(\|\theta_0^n\|_{B^{-1}_{\infty,1}}
+\int_0^t\|u^n\cdot\nabla\theta^n\|_{B^{-1}_{\infty,1}}\,d\tau\biggr).
\end{equation}
In order to bound the source term, one may use
the following \emph{Bony's  decomposition}: 
\begin{equation}\label{eq:bony}
u^n\cdot\nabla\theta^n=\div
R(u^n,\theta^n)+\sum_{j=1}^2\Bigl(T_{\d_j\theta^n}u^n_j+ T_{u^n_j}\d_j\theta^n\Bigr).
\end{equation}
In the above formula, $T$ (resp. $R$) stands for the paraproduct
(resp. remainder) operator  defined by 
\begin{equation}\label{eq:defpara}
T_fg:=\sum_{q\geq1} S_{q-1}f\Delta_qg\quad\biggl(\hbox{resp. }\
R(f,g):=\sum_{q\geq-1} \Delta_q
f\,\tilde\dq g\biggr)
\end{equation}
with $S_p:=\sum_{p'\leq p-1}\Delta_{p'}$ and
$\tilde\Delta_p:=\Delta_{p-1}+\Delta_p+\Delta_{p+1},$
and we use the fact that, owing to  $\div u^n=0,$
we have
$$
\sum_{j=1}^2R(u^n_j,\d_j\theta^n)=\div R(u^n,\theta^n).
$$
For the remainder term $R,$ it is standard (see e.g. \cite{BCD}) that 
\begin{equation}\label{eq:remainder}
\|R(u^n,\theta^n)\|_{B^{1}_{\infty,\infty}} \leq
C\|\theta^n\|_{B^{0}_{\infty,\infty}}\|u^n\|_{B^1_{\infty,\infty}}.
\end{equation} Now, because  $\Delta u^n=\nabla^\bot \omega^n$ with
$\nabla^\bot:=(-\d_2,\d_1),$
one may decompose $u^n$ into 
$$
u^n=\Delta_{-1}u^n-\sum_{q\geq0}\nabla^\bot(-\Delta)^{-1}\Delta_q\omega^n.
$$
Using Bernstein inequalities and the fact that operator
 $\nabla^\bot(-\Delta)^{-1}$ is homogeneous of degree $-1,$ we eventually get
\begin{equation}\label{eq:u}
\|u^n\|_{B^1_{\infty,\infty}}\leq C\bigl(\|u^n\|_{L^\infty}+
\|\omega^n\|_{L^\infty}\bigr).
\end{equation}
As operator $\div$ maps
$B^1_{\infty,\infty}$ in $B^0_{\infty,\infty},$  and
as $B^0_{\infty,\infty}\hookrightarrow B^{-1}_{\infty,1}$ and
$H^1\hookrightarrow B^{0}_{\infty,\infty},$  
we thus get from \eqref{eq:remainder} and \eqref{eq:u},
\begin{equation}\label{eq:para1} \|\div R(\theta^n,u^n)\|_{B^{-1}_{\infty,1}}
\leq C\|\theta^n\|_{H^1}\bigl(\|u^n\|_{L^\infty}+\|\omega^n\|_{L^\infty}\bigr).
\end{equation}
Next, making use of continuity properties for the paraproduct operator (see
e.g. \cite{BCD}), we discover that $$
\|T_{\d_j\theta^n}u^n_j\|_{B^{-1}_{\infty,1}}+
\|T_{u^n_j}\d_j\theta^n\|_{B^{-1}_{\infty,1}}\leq
C\|u^n_j\|_{L^\infty}\|\d_j\theta^n\|_{B^{-1}_{\infty,1}}\quad\hbox{for}\quad
j=1,2.
$$
Plugging this latter inequality and \eqref{eq:para1} in \eqref{eq:bony},
we get 
\begin{equation}\label{eq:para2}\|u^n\cdot\nabla\theta^n\|_{B^{-1}_{\infty,1}}
\leq 
C\Bigl(\bigl(\|u^n\|_{L^\infty}+\|\omega^n\|_{L^\infty}\bigr)\|\theta^n\|_{H^1} +\|u^n\|_{L^\infty}\|\theta^n\|_{B^{0}_{\infty,1}}\Bigr).
\end{equation}
In order to conclude, one may use the following
two inequalities the proof of which has been  postponed in the appendix:
\begin{eqnarray}\label{eq:interpo1}
&\|u^n\|_{L^\infty}\leq C\|u^n\|_{L^2}^{\frac 12}
\|\omega^n\|_{L^\infty}^{\frac 12},\\\label{eq:interpo2}
&\|\theta^n\|_{B^{0}_{\infty,1}}
\leq C\|\theta^n\|_{L^2}^{\frac12}\|\theta^n\|_{B^1_{\infty,1}}^{\frac12}.
\end{eqnarray}
Inserting \eqref{eq:interpo1} and \eqref{eq:interpo2} in \eqref{eq:para2}
then using Young inequality, we get for all $\varepsilon>0,$
$$\displaylines{
\int_0^t\|u^n\cdot\nabla\theta^n\|_{B^{-1}_{\infty,1}}\,d\tau
\leq
C\biggl(\int_0^t\|\theta^n\|_{H^1}\bigl(\|u^n\|_{L^2}
+\|\omega^n\|_{L^\infty}\bigr)\,d\tau
\hfill\cr\hfill+\frac{1+\kappa
t}{\varepsilon\kappa}\int_0^t\|u^n\|_{L^2}\|\omega^n\|_{L^\infty}
\|\theta^n\|_{L^2}\,d\tau+\frac{\varepsilon\kappa}{1+\kappa t}
\int_0^t\|\theta^n\|_{B^{1}_{\infty,1}}\,d\tau\biggr). }
$$
Taking 
 $\varepsilon$ sufficiently small and coming back to \eqref{eq:chaleur}, 
we end up with
 $$\displaylines{
\Theta^n(t)\leq C(1+\kappa t)\biggl(\Theta^n_0
+\int_0^t\|\theta^n\|_{H^1}\|u^n\|_{L^2}\,d\tau\hfill\cr\hfill
+\int_0^t\bigl(\|\theta^n\|_{H^{1}}+(\kappa^{-1}+t)\|u^n\|_{L^2}\|\theta^n\|_{L^2}\bigl)
\|\omega^n\|_{L^\infty}\,d\tau\biggr)}
$$
where $\Theta^n(t):=\sup_{\alpha\in[1,\infty]} 
\kappa^{\frac1\alpha}\|\theta^n\|_{L_t^\alpha(B^{-1+\frac2\alpha}_{\infty,1})}$
and $\Theta^n_0:=\|\theta_0^n\|_{B^{-1}_{\infty,1}}.$
\smallbreak
On the one hand, the above inequality rewrites
\begin{equation}\label{eq:borne0}
\Theta^n(t)\leq f^n(t)+ (1+\kappa t)^2
\int_0^t g^n(\tau)\|\omega^n(\tau)\|_{L^\infty}\,d\tau
\end{equation}
with $\ \left\{\begin{array}{lll}
f^n(t)&=& C(1+\kappa t)\biggl(\Theta_0^n
+\Int_0^t\|\theta^n\|_{H^1}\|u^n\|_{L^2}\,d\tau\biggr),\\[1.5ex]
g^n(t)&=&C\bigl(\|\theta^n(t)\|_{H^1}
+\kappa^{-1}\|u^n(t)\|_{L^2}\|\theta^n(t)\|_{L^2}\bigr).
\end{array}\right.$
\medbreak
On the other hand, according to  \eqref{eq:tourbillonLp}
and  as $B^0_{\infty,1}\hookrightarrow L^\infty,$ we have
\begin{equation}\label{eq:omega}
\|\omega^n(t)\|_{L^\infty}\leq \|\omega_0^n\|_{L^\infty}
+C\kappa^{-1}\Theta^n(t).
\end{equation}
Inserting the above inequality in \eqref{eq:borne0}
and making use of Gronwall lemma thus  yields
\begin{equation}\label{eq:borne1}
\Theta^n(t)\leq \biggl(f^n(t)+(1+\kappa t)^2\|\omega_0^n\|_{L^\infty}
\int_0^tg^n(\tau)\,d\tau\biggr)
e^{C\kappa^{-1}(1+\kappa t)^2\int_0^tg^n(\tau)\,d\tau}.
\end{equation}
Obviously, \eqref{eq:thetaL2} and $\eqref{eq:uL2}$ imply that 
$(u^n)_{n\in\N}$ is bounded in $L^\infty_{loc}(\R_+;L^2)$
and that  $(\theta^n)_{n\in\N}$ is bounded in
$L^\infty(\R_+;L^2)\cap L^2_{loc}(\R_+;H^1).$
Therefore the right-hand side of \eqref{eq:borne1}
may be  bounded independently of $n.$ 
This provides a uniform bound for  $\theta^n$ in the space
$L^1_{loc}(\R_+;B^1_{\infty,1})\cap  L^\infty_{loc}(\R_+;B^{-1}_{\infty,1}).$
Next, coming back to  \eqref{eq:omega} yields
a bound  for $(\omega^n)_{n\in\N}$ in $L^\infty_{loc}(\R_+;L^\infty).$

\subsubsection*{Fourth step.} In order to show that 
$(\theta^n,u^n)_{n\in\N}$ converges (up to extraction), 
a boundedness information over 
 $(\d_t\theta^n,\d_tu^n)$ is needed. 

As for the temperature, because
$$
\d_t\theta^n=\kappa\Delta\theta^n-u^n\cdot\nabla\theta^n,
$$
the previous steps  imply that $(\d_t\theta^n)_{n\in\N}$
 is bounded in $L^2_{loc}(\R_+;H^{-1}).$ 

We claim that
$(\d_tu^n)_{n\in\N}$ is bounded in
$L^\infty_{loc}(\R_+;L^2).$
Indeed, applying the Leray projector
 ${\mathcal P}$ over divergence free vector-fields to the velocity equation 
yields
$$
\d_tu^n=-{\mathcal P}(\theta^n e_2-u^n\cdot\nabla u^n).
$$
Since  $(\theta^n)_{n\in\N}$ is bounded in $L^\infty(\R_+;L^2),$ so is 
 ${\mathcal P}(\theta^n e_2).$  Next, as  $(\omega^n)_{n\in\N}$ 
is bounded in $L^\infty_{loc}(\R_+;L^r),$   so is 
$(\nabla u^n)_{n\in\N}$ according to proposition \ref{p:biot-savart}. Finally,
the previous results imply that sequence $(u^n)_{n\in\N}$ is bounded in
$L^\infty_{loc}(\R_+;L^2\cap L^\infty),$ thus in $L^\infty_{loc}(\R_+;L^s)$
with  $s=2r/(r-2).$ Thanks to H\"older inequality, one can thus conclude that
  $(u^n\cdot\nabla u^n)_{n\in\N}$ is bounded in $L^\infty_{loc}(\R_+;L^2).$

\subsubsection*{Fifth step.} Passing to the limit.

According to the previous steps, we have
\begin{itemize}
\item $(\theta^n)_{n\in\N}$ is bounded in 
$L^\infty_{loc}(\R_+;L^2\cap B^1_{\infty,1})
\cap L^2_{loc}(\R_+;H^1)\cap L^1_{loc}(\R_+;B^1_{\infty,1}),$
\item $(\d_t\theta^n)_{n\in\N}$ is bounded in $L^2_{loc}(\R_+;H^{-1}),$
\item $(u^n)_{n\in\N}$ and $(\d_tu^n)_{n\in\N}$ are bounded in
$L^\infty_{loc}(\R_+;L^2),$ 
\item $(\omega^n)_{n\in\N}$ is bounded in  $L^\infty_{loc}(\R_+;L^r\cap
L^\infty).$ 
\end{itemize}
Because  $H^{-1}$ is (locally) compactly embedded in $L^2$
the classical Aubin-Lions argument  (see e.g.  \cite{Aubin})
ensures  that, up to extraction, 
sequence  $(\theta^n,u^n)_{n\in\N}$ strongly converges in 
$L^\infty_{loc}(\R_+;H^{-1}_{loc})$
to some function  $(\theta,u)$ so that 
$$\displaylines{
\theta\in L^\infty_{loc}(\R_+;L^2\cap B^1_{\infty,1})
\cap L^2_{loc}(\R_+;H^1)\cap L^1_{loc}(\R_+;B^1_{\infty,1}),\cr
u\in C^{0,1}_{loc}(\R_+;L^2)\!\!\quad\hbox{and}\!\!\quad
\omega\in L^\infty_{loc}(\R_+;L^r\cap L^\infty).}
$$
Now, interpolating with the uniform bounds stated in the previous
steps, it is  easy to pass to the limit in $(B_{\kappa,0}).$  
Finally, from standard properties for the heat equation  (see e.g. \cite{CH})
we get in addition  $\theta\in\cC(\R_+;L^2\cap B^{-1}_{\infty,1}).$
This completes the proof of existence.

\subsubsection*{Sixth step.} In order to show the uniqueness part
of our statement, we shall use the Yudovich argument \cite{Yudo}
revisited by P. G\'erard in \cite{PG}.

Let  $(\theta_1,u_1,\Pi_1)$ and  
$(\theta_2,u_2,\Pi_2)$ satisfy \eqref{eq:regularity} and
 $(B_{\kappa,0})$ with
the same data. 
Denote $\dt:=\theta_2-\theta_1,$ $\du:=u_2-u_1$ 
and  $\dPi:=\Pi_2-\Pi_1.$
Because 
$$
\d_t\du+u_2\cdot\nabla\du+\nabla\dPi=-\du\cdot\nabla u_1+\dt\,e_2,
$$ 
a standard energy method combined with  H\"older inequality yields
for all  $p\in[2,\infty[$   $$
\frac12\frac d{dt}\|\du\|_{L^2}^2
\leq \|\nabla u_1\|_{L^p}\|\du\|_{L^{2p'}}^2+\|\dt\|_{L^2}\|\du\|_{L^2}\quad\hbox{with }\
p':=\frac{p}{p-1}\cdotp
$$
This inequality rewrites
\begin{equation}\label{eq:du}
\frac12\frac d{dt}\|\du\|_{L^2}^2
\leq p\|\nabla
u_1\|_{L}\|\du\|_{L^\infty}^{\frac2p}
\|\du\|_{L^{2}}^{\frac2{p'}}+\|\dt\|_{L^2}\|\du\|_{L^2}
\end{equation}
 with
$$
\|\nabla u_1\|_{L}:=\sup_{r\leq p<\infty}\frac{\|\nabla u_1\|_{L^p}}p\cdotp
$$
Let us point out that, by virtue of Proposition \ref{p:biot-savart},
as $\omega_1\in L^\infty_{loc}(\R_+;L^r\cap L^\infty)$
the term $\|\nabla u_1(t)\|_{L}$ is locally bounded.
Of course, combining the fact that
 $u_i\in L^\infty_{loc}(\R_+;L^2)$ and  $\omega_i\in
L^\infty_{loc}(\R_+;L^\infty)$ for $i=1,2,$  implies that   $\du\in
L^\infty_{loc}(\R_+;L^\infty).$ 
\smallbreak
Next, we notice that  $\dt$ satisfies
$$
\d_t\dt-\kappa\Delta\dt=-u_2\cdot\nabla\dt-\du\cdot\nabla\theta_1,\qquad
\d_t\dt_{|t=0}=0.
$$ 
Our regularity assumptions over the solutions ensure that the right-hand
side belongs to  $L^2_{loc}(\R_+;L^2)$.
  Hence, according to a standard maximal regularity result for the heat equation, 
  we deduce that  $\d_t\delta\theta\in L^2_{loc}(\R_+;L^2).$
   Hence, using an energy method yields
\begin{equation}\label{eq:dt}\frac12
\frac d{dt}\|\dt\|_{L^2}^2\leq
\|\nabla\theta_1\|_{L^\infty}\|\dt\|_{L^2}\|\du\|_{L^2}.
\end{equation}
Let $\e$ be a small parameter  (bound to tend to $0$). 
Denote
$$
X_\e(t):=\sqrt{\|\dt(t)\|_{L^2}^2+\|\du(t)\|_{L^2}^2+\e^2}.
$$
Putting inequalities \eqref{eq:du} and  \eqref{eq:dt} together gives
$$
\frac d{dt} X_\e\leq p\|\nabla u_1\|_{L}\|\du\|_{L^\infty}^{\frac2p}
X_\e^{1-\frac2{p}}+\frac12
(1+\|\nabla\theta_1\|_{L^\infty})X_\e.
$$
Let $\gamma(t):=\frac12(1+\|\nabla\theta_1(t)\|_{L^\infty}).$
The assumptions over $\theta_1$ ensure that
function  
$\gamma$ is in $L^1_{loc}(\R_+).$ 
Therefore, setting
$Y_\e:=e^{-\int_0^t\gamma(\tau)\,d\tau}X_\e,$
the previous inequality rewrites
$$
\frac2p Y_\e^{\frac2p-1}\frac d{dt} Y_\e\leq 2\|\nabla
u_1\|_{L}\|\du\|_{L^\infty}^{\frac2p}
e^{-\frac2p\int_0^t\gamma(\tau)\,d\tau}. $$
Performing a time integration yields
$$
Y_\e(t)\leq\biggl(\e^{\frac2p}
+2\int_0^t\|\nabla
u_1\|_{L}\|\du\|_{L^\infty}^{\frac2p}\,d\tau\biggr)^{\frac p2}.
$$
Having $\e$ tend to $0,$ we end up with
\begin{equation}\label{eq:dt1}
\|\dt(t)\|_{L^2}^2+\|\du(t)\|_{L^2}^2
\leq \|\du\|_{L_t^\infty(L^\infty)}^2 \biggl(2\int_0^t\|\nabla
u_1\|_{L}\,d\tau\biggr)^p\quad\hbox{for all }\ t\in\R^+.
\end{equation}
As explained above, the term 
$\|\nabla u_1(t)\|_{L}$ is locally bounded. Hence
one may find  a positive time $T$ 
so that   $\int_0^T\|\nabla
u_1\|_{L}\,d\tau<\frac12.$ 
Letting  $p$ tend to infinity in \eqref{eq:dt1} thus entails
that  $(\dt,\du)\equiv0$ on $[0,T].$
Because $\dt$ and $\du$ are continuous in time with values in $L^2,$
it is now easy to conclude that $(\dt,\du)\equiv0$ on $\R^+,$
by means of a standard connectivity argument.


\section{A global result for infinite energy initial velocity}

In dimension two, the assumption that 
$u_0$ is in $L^2$ is somewhat restrictive since it entails that
the vorticity $\omega_0$ has $0$ average over $\R^2.$
This in particular precludes our considering vortex patches like structures
or, more generally, data with compactly supported
nonnegative vorticity.
The present section aims at generalizing our study to initial velocity fields
with  (possibly) infinite energy. 
The functional setting we shall introduce below is borrowed 
from Chemin's in \cite{C}. 

Let us first notice that whenever $g$ is a radial $\cC_c^\infty$ function supported
away from the origin  then the smooth vector field $\sigma$ defined 
by  \begin{equation}\label{eq:sigma}
\sigma(x)=\frac{x^\bot}{|x|^2}\int_0^{|x|} rg(r)\,dr\end{equation}
is a \emph{stationary} solution to the two-dimensional incompressible Euler
equations, and has vorticity $\omega_\sigma: x\mapsto g(|x|).$

For $m\in\R,$ we then define $E_m$ as the set of all divergence-free $L^2$
perturbations of a velocity field $\sigma$ satisfying \eqref{eq:sigma} and
\begin{equation}\label{eq:sigma1}
\int_{\R^2} g(|x|)\,dx=m.
\end{equation}
Showing that the definition of $E_m$ depends 
only on $m$ is left to the reader (it is only a matter of using Fourier
variables).  \smallbreak
The rest of this section is devoted to the proof of the following
generalization of Theorem~\ref{th:main1}.
\begin{theorem}\label{th:main2}
Let  $\theta_0\in L^2\cap B^{-1}_{\infty,1}$ and  $u_0\in E_m$ for some $m\in\R.$
Assume in addition that the initial vorticity $\omega_0$ belongs to 
$L^r\cap L^\infty$ for some  $r\geq2.$
 Then System  $(B_{\kappa,0})$ admits a unique
global solution  $(\theta,u)$ such that
\begin{equation}\label{eq:regularity2}
\begin{array}{c}
\theta\in\cC(\R_+;L^2\cap
B^1_{\infty,1}) \cap L^2_{loc}(\R_+;H^1)\cap L^1_{loc}(\R_+;B^1_{\infty,1}),\\[1.5ex]
u\in\cC^{0,1}_{loc}(\R_+;E_m) \ \hbox{ and }\ \omega\in
L^\infty_{loc}(\R_+;L^r\cap L^\infty).\end{array}\end{equation}  \end{theorem} 
 \begin{p} 
 As it is very similar to that of Theorem \ref{th:main1}, 
 we just sketch the proof and
 point out what has to be changed. 
 
 Throughout we fix a stationary vector-field  $\sigma$ 
 satisfying  \eqref{eq:sigma} and \eqref{eq:sigma1}. Setting
 $u=v+\sigma,$ System $(B_{\kappa,0})$ rewrites
 \begin{equation}\label{eq:modifiee}
 \begin{cases} \partial_t\theta+(v+\sigma)\cdot\nabla
\theta-\kappa\Delta \theta=0\\ \partial_t v+(v+\sigma)\cdot\nabla v+v\cdot\nabla\sigma
+\nabla\Pi=\theta\,
e_2\\ \div v=0. \end{cases} 
  \end{equation}
 As $\div\sigma=\div v=0,$ the energy estimates for $\theta$ remain the same. 
As for the velocity field, having the new term $v\cdot\nabla\sigma$ in the equation
implies that
 \begin{equation}\label{eq:vL2}
 \|v(t)\|_{L^2}\leq e^{t\|\nabla \sigma\|_{L^\infty}}
 \|v_0\|_{L^2}+\biggl(\frac{e^{t\|\nabla\sigma\|_{L^\infty}}-1}{\|\nabla\sigma\|_{L^\infty}}
 \biggr)\|\theta_0\|_{L^2}.
 \end{equation}
 Now, the vorticity $\omega_v$ associated to  $v$ satisfies
 $$
 \d_t\omega_v+(v+\sigma)\cdot\nabla\omega_v+v\cdot\nabla\omega_\sigma=\d_1\theta.
 $$
 Hence for all  $p\in[r,\infty],$
 $$
 \|\omega_v(t)\|_{L^p}\leq\|\omega_v(0)\|_{L^p}
 +\int_0^t\|\d_1\theta\|_{L^p}\,d\tau
 +\int_0^t\|v\|_{L^p}\|\nabla\omega_\sigma\|_{L^\infty}\,d\tau.
 $$
 Splitting $v$ into  
 $$
 v=\Delta_{-1}v-\sum_{q\in\N}\nabla^\bot(-\Delta)^{-1}\Delta_q\omega_v
 $$
and using Bernstein inequality,  we readily get
 $$
 \|v\|_{L^p}\leq C\bigl(\|v\|_{L^2}+\|\omega_v\|_{L^p}\bigr).
 $$
 Therefore, as in the proof of theorem
 \eqref{th:main1}, in order to bound $\omega_v$ in 
$L^\infty_{loc}(\R_+;L^r\cap L^\infty),$ it suffices
to get  a  bound for
 $\d_1\theta$ in $L^1_{loc}(\R_+;L^\infty).$
  This may be achieved by bounding  $\d_1\theta$ in
$L^1_{loc}(\R_+;B^0_{\infty,1}),$ given  that
 $$
 \d_t\theta-\kappa\Delta\theta=-v\cdot\nabla\theta-\sigma\cdot\nabla\theta.
 $$
Arguing as in \eqref{eq:chaleur} reduces the problem to getting
 an appropriate  bound for the new term 
$\sigma\cdot\nabla\theta$ in  $L^1_{loc}(\R_+;B^{-1}_{\infty,1}).$  For this
purpose, one may use again Bony's decomposition, the fact that $\div\sigma=0$
and classical continuity properties for the paraproduct and remainder
operators. One ends up for instance with:  $$
 \|\sigma\cdot\nabla\theta\|_{B^{-1}_{\infty,1}}\leq C\|\sigma\|_{B^\e_{\infty,\infty}}
 \|\theta\|_{L^\infty}.
 $$
 Combining \eqref{eq:interpo1} and Young inequality, it is now easy to 
 get an inequality similar to  \eqref{eq:borne0}, and thus a bound for $\theta$
in  $L^1_{loc}(\R_+;B^1_{\infty,1})\cap
L^\infty_{loc}(\R_+;B^{-1}_{\infty,1}).$  \smallbreak
 In order to prove the uniqueness, it is fundamental to notice that
 if  $(\theta_1,u_1)$ and  $(\theta_2,u_2)$ both solve 
  $(B_{\kappa,0})$ with the same data, and satisfy 
  \eqref{eq:regularity2} \emph{with the same $m$}
(an assumption which is not restrictive since we know that $u_1$
and $u_2$ coincide initially) then one may 
  write  $u_1=\sigma+v_1$ and  $u_2=\sigma+v_2$ for some 
 stationary vector-field $\sigma$ satisfying  \eqref{eq:sigma},\eqref{eq:sigma1}
and  $v_1,v_2$ in $L^\infty_{loc}(\R_+;L^2).$
 
 Taking advantage of Equation $\eqref{eq:modifiee}_2,$ is is obvious that  
 $\d_tv_1$ and  $\d_tv_2$ are in   $L^\infty_{loc}(\R_+;L^2).$
 Now, we notice that $(\dv,\dt):=(v_2-v_1,\theta_2-\theta_1)$ satisfy 
 $$
 \left\{\begin{array}{l}
 \d_t\dt+u_2\cdot\nabla\dt-\kappa\Delta\dt=-\dv\cdot\nabla\theta_1,\\[1ex]
\d_t\dv+u_2\cdot\nabla\dv+\nabla\dPi=-\dv\cdot\nabla u_1+\dt\,e_2-\dv\cdot\nabla\sigma.
\end{array}\right.
$$ 
Up to the additional term  $-\dv\cdot\nabla\sigma$
which may be bounded as follows:
$$
\|\dv\cdot\nabla\sigma\|_{L^2}\leq\|\dv\|_{L^2}\|\nabla\sigma\|_{L^\infty},
$$ 
the energy bounds for the above system 
are the same as in the case $\sigma=0.$
 Hence, from argument similar to those used in the previous section, 
 it is easy to conclude the proof of uniqueness. 
 The details are left to the reader.
 \end{p}


\section{Further results and concluding remarks}

In this concluding section, we list a few extensions
which may be obtained by straightforward generalizations
of our method.
\subsection{Remarks concerning the Boussinesq system}
Let us stress that the key to the proof of Theorems \ref{th:main1} and \ref{th:main2}
is that, on the one hand,  the solution does not develop singularities as long
as  $$
\int_0^T\|\nabla\theta\|_{L^\infty}\,dt<\infty,
$$
and that, on the other hand, under quite weak assumptions over the initial
data,  the above integral remains finite for all $T<\infty.$
\smallbreak
In fact, a quick revisitation of our proof shows that 
if one assumes in addition that 
   $\omega_0\in C^\epsilon$ and  $\theta_0\in C^{-1+\epsilon}$
(with $C^{-1+\epsilon}:=B^{-1+\epsilon}_{\infty,\infty}$)
 for some  $\e\in]0,1[$ then both
 $\nabla\theta$ and $\nabla u$ are in $L^1_{loc}(\R_+;L^\infty(\R^2))$
so that  the additional H\"older regularity is conserved during the evolution.
We believe that, more generally, our study 
opens a way to  investigate  vortex patches structures 
(or striated regularity) for the Boussinesq system with $\kappa>0$ and $\nu=0.$
 \medbreak
Let us also emphasize that 
if, in addition to the hypotheses of Theorem \ref{th:main1},
we have $u_0\in B^1_{\infty,1}$ then the corresponding solution $(\theta,u)$
also satisfies
 $$ u\in\cC(\R_+;B^1_{\infty,1}).
$$
 Indeed, according to a result by M. Vishik in
\cite{Vishik} concerning  the transport equation, one can propagate the
$B^0_{\infty,1}$ regularity over the vorticity $\omega$ provided  $\d_1\theta$
is in   $L^1_{loc}(\R_+; B^0_{\infty,1})$  
and  there exists some universal constant $C$ such that
\begin{equation}\label{eq:besov}
\|\omega(t)\|_{B^0_{\infty,1}}\leq C
\biggl(1+\int_0^t\|\nabla u\|_{L^\infty}\biggr)
\biggl(\|\omega_0\|_{B^0_{\infty,1}}+\int_0^t\|\d_1\theta\|_{B^0_{\infty,1}}
\biggr).
\end{equation}
  Now, under the sole  assumptions   of Theorem \ref{th:main1},
one may bound $\d_1\theta$ in  
$L^1_{loc}(\R_+; B^0_{\infty,1})$ by means of the norms of the data.
Because, owing to $B^0_{\infty,1}\hookrightarrow L^\infty$
and \eqref{eq:interpo1}, one may write
$$
\|\nabla u\|_{L^\infty}\leq
C\bigl(\|u\|_{L^2}+\|\omega\|_{B^0_{\infty,1}}\bigr),
$$
Inequality \eqref{eq:besov} combined with Gronwall lemma
ensures  the conservation of
the additional $B^0_{\infty,1}$ regularity for the vorticity
(and thus of the $B^1_{\infty,1}$ regularity for the velocity).
 This argument provides another proof of Hmidi and 
Keraani's result in  \cite{HK2} under somewhat weaker assumptions over
$\theta_0$ (there having   $\theta_0$ in (a subspace of) $L^\infty$ was
needed).

\subsection{The B\'enard system}

Our method may also be adapted with almost no change
to the study of the following \emph{B\'enard system}: 
\begin{equation}\label{eq:benard}
\begin{cases}
\partial_t\theta+u\cdot\nabla\theta-\kappa\Delta\theta=u_2\\
\partial_t u+u\cdot\nabla u+\nabla p=\theta\, e_2\\
(\theta, u)|_{t=0}=(\theta_0, u_0),
\end{cases}
\end{equation}
which describes convective motions in a heated 
 two-dimensional inviscid incompressible
fluid under thermal effects (see e.g. \cite{AP}, Chap. 6). We get
\begin{theorem}\label{th:main3}
For all data  $(\theta_0,u_0)$
with  $\theta_0\in L^2\cap B^{-1}_{\infty, 1}$ and  
 $u_0\in L^2$ satisfying  $\div u_0=0$ and   $\omega_0\in L^r\cap L^\infty$ for some
 $r\in[2,\infty[,$  System $\eqref{eq:benard}$ has a unique global solution
   $(\theta,u)$ such that
  \begin{equation}\label{eq:regularity3}
  \begin{array}{c}\theta\in\cC(\R_+;L^2\cap
B^1_{\infty,1}) \cap L^2_{loc}(\R_+;H^1)\cap L^1_{loc}(\R_+;B^1_{\infty,1}),\\[1.5ex]
\quad u\in\cC^{0,1}_{loc}(\R_+;L^2)\  \hbox{ and }\ \omega\in
L^\infty_{loc}(\R_+;L^r\cap L^\infty).\end{array}\end{equation}\end{theorem} 
\begin{p}
We just briefly indicate what has to be changed compared 
to the proof of Theorem~\ref{th:main1}.
Owing to the new term $u_2$ in the equation for the temperature, 
the energy estimates read
\begin{eqnarray}\label{eq:benard1}\frac{1}{2}\frac{d}{dt}\|\theta\|_{L^2}^2
+\kappa\|\nabla \theta\|_{L^2}^2=\int\theta\: u_2\, dx,\\
\label{eq:benard2}
\frac 12 \frac{d}{dt}\|u\|_{L^2}^2=\int\theta \: u_2\,dx.
\end{eqnarray}
 Adding up inequalities \eqref{eq:benard1} and \eqref{eq:benard2} yields
$$
\Frac{1}{2}\frac{d}{dt}\|(\theta,u)(t)\|_{L^2}^2
+\kappa\|\nabla \theta\|_{L^2}^2=2\Int\theta\:u_2\,dx,\\[2ex]
\leq\|(\theta,u)\|_{L^2}^2.
$$
Thanks to the Gronwall inequality, we thus infer that 
$$\|(\theta,u)(t)\|_{L^2}^2+2\kappa\int_0^t \|\nabla
\theta(\tau)\|_{L^2}^2\,d\tau \leq \|(\theta_0,u_0)\|_{L^2}^2\:e^{2t}.$$
The rest of the proof of Theorem \ref{th:main3} follows the lines
of that of Theorem \ref{th:main1}, once it has been noticed
that the computations  leading to Inequality  \eqref{eq:chaleur} (see the
appendix) also yield  $$\bigg\|\int_0^t e^{(t-s)\kappa\Delta} u_2(s)
ds\bigg\|_{L_T^1(B^{1}_{\infty,1})}\leq   C\int_0^T \|u\|_{L^\infty}\,dt.
 $$
 Note also that having the new (lower order) term $u_2$ 
 in Equation $\eqref{eq:benard}_1$ is harmless 
 for proving uniqueness.
\end{p}


\section*{Appendix}

Here we prove a few inequalities which have been used
throughout  the paper. 
\smallbreak\noindent{\it Proof of Inequality $\eqref{eq:chaleur}$:}
  Assume that  $\theta$ satisfies 
  $$
  \d_t\theta-\kappa\Delta\theta=f,\qquad \theta_{|t=0}=\theta_0.
  $$
  Then applying the dyadic operator $\Delta_q$ to the above equality
  yields
    $$
  \d_t\dq\theta-\kappa\dq\Delta\theta=\dq f\quad\hbox{for all }\ q\geq-1.
  $$
  {}From the maximum principle, we readily get
  $$
\|\Delta_{-1}\theta(t)\|_{L^\infty}\leq \|\Delta_{-1}\theta_0\|_{L^\infty}
+\int_0^t\|\Delta_{-1}f(\tau)\|_{L^\infty}\,d\tau
$$
whence for all  $\alpha\in[1,\infty]$ and  $t>0,$
\begin{equation}\label{eq:heat1}
 \|\dq\theta\|_{L^\alpha([0,t];L^\infty)}\leq
Ct^{\frac1\alpha}\Bigl(\|\Delta_{-1}\theta_0\|_{L^\infty}+\|\Delta_{-1}f\|_{L^1([0,t];L^\infty)}\Bigr). \end{equation}
Next, for bounding the high frequency blocks $\Delta_q\theta$ with  $q\geq0,$ 
one may write  
 \begin{equation}\label{eq:heat2}
  \Delta_q\theta(t)=e^{\kappa t\Delta}\dq\theta_0
  +\int_0^te^{\kappa(t-\tau)\Delta}\dq f(\tau)\,d\tau
  \end{equation}
 where $(e^{\lambda\Delta})_{\lambda>0}$ stands 
  for the heat semi-group, and take advantage of the 
  following inequality stated by J.-Y. Chemin in \cite{CH}:
there exists two positive constants $c$ and $C$ 
such that  
  \begin{equation}\label{eq:chemin}
  \|e^{\lambda\Delta}\dq g\|_{L^\infty}\leq
   Ce^{-c\lambda2^{2q}}\|\dq g\|_{L^\infty}\ \hbox{ for all }\ \lambda>0\
\hbox{ and }\ q\geq0.    \end{equation}
  {}From \eqref{eq:heat2} and \eqref{eq:chemin}, we get
  $$
  \|\Delta_q\theta(t)\|_{L^\infty}\leq C\biggl(e^{-c\kappa2^{2q}t}\|\Delta_q\theta_0\|_{L^\infty}
  +\int_0^te^{-c\kappa2^{2q}(t-\tau)}\|\Delta_qf(\tau)\|_{L^\infty}\,d\tau\biggr).
  $$
Therefore, for all  $\alpha\in[1,\infty],$ $q\geq0$ and  $t>0,$ 
 $$
\kappa^{\frac1\alpha}2^{(\frac2\alpha-1)q}\|\dq\theta\|_{L^\alpha([0,t];L^\infty)}
\leq  C2^{-q}\Bigl(\|\dq\theta_0\|_{L^\infty}+\|\dq
f\|_{L^1([0,t];L^\infty)}\Bigr).  $$ Summing on  $q\geq0$ and using 
\eqref{eq:heat1}, it is now easy to complete the proof of Inequality 
\eqref{eq:chaleur}.     \hfill\rule{2.1mm}{2.1mm} 
 \smallbreak\noindent{\it Proof of Inequalities $\eqref{eq:interpo1}$
 and  $\eqref{eq:interpo2}$:}
 For proving the first inequality, let us consider
 a $L^2$ divergence free vector-field $u$  with bounded vorticity  $\omega.$ 
 As  $u$ is in  $L^2,$ one may write
$$
u=\sum_{q\in\Z}\dot\Delta_qu\quad\hbox{with}\quad
\dot\Delta_q:=\varphi(2^{-q}D).
$$
Let $N$ be an integer parameter to be chosen hereafter. 
Given that 
$u=-\nabla^{\perp}(-\Delta)^{-1}\omega$ and using the Bernstein
inequalities, we have
 $$\|u\|_{L^\infty}\leq \sum\limits_{q\leq N}\|\dot\Delta_q 
u\|_{L^\infty}+\sum\limits_{q>N}\|\dot \Delta_q u\|_{L^\infty}\leq  
C2^N\|u\|_{L^2}+C\sum\limits_{q>N}2^{-q}
\|\dot\Delta_q\omega\|_{L^\infty}.$$
Therefore, 
 $$\|u\|_{L^\infty}\leq 
C2^N\|u\|_{L^2}+C2^{-N}\|\omega\|_{L^\infty}.$$
Taking  $N$ so that  $2^N\|u\|_{L^2}
\approx2^{-N}\|\omega\|_{L^\infty}$, 
we get the desired inequality.

Proving Inequality \eqref{eq:interpo2} relies on the
similar decomposition into low and high frequencies. The 
details are left to the reader. 
  \hfill\rule{2.1mm}{2.1mm}


\end{document}